\documentclass[review]{elsarticle}

\usepackage{hyperref}

\usepackage{amsmath}
\usepackage{amsfonts}
\usepackage{siunitx}
\usepackage{xcolor}
\usepackage{color,soul}

\begin{document}

\begin{frontmatter}

\title{Reinforcement Learning Versus Model Predictive Control on Greenhouse Climate Control}

\author{Bernardo Morcego\fnref{myfootnote1}}
\author{Wenjie Yin\fnref{myfootnote1}}
\author{Sjoerd Boersma\fnref{myfootnote2}}
\author{Eldert van Henten\fnref{myfootnote3}}
\author{Vicen\c c Puig\fnref{myfootnote1}}
\author{Congcong Sun\fnref{myfootnote3}\texorpdfstring{\corref{cor1}}{}}

\fntext[myfootnote1]{Automatic Control Group (CS2AC), Universitat Politècnica de Catalunya, Spain}
\fntext[myfootnote2]{Biometris Group, Wageningen University, 6700 AA Wageningen, The Netherlands}
\fntext[myfootnote3]{Agricultural Biosystems Engineering Group, Wageningen University, 6700 AA Wageningen, The Netherlands}
\cortext[cor1]{Corresponding email: congcong.sun@wur.nl}

\begin{abstract}
Greenhouse is an important protected horticulture system for feeding the world with enough fresh food. However, to maintain an ideal growing climate in a greenhouse requires resources and operational costs. In order to achieve economical and sustainable crop growth, efficient climate control of greenhouse production becomes essential. Model Predictive Control (MPC) is the most commonly used approach in the scientific literature for greenhouse climate control. However, with the developments of sensing and computing techniques, reinforcement learning (RL) is getting increasing attention recently. With each control method having its own way to state the control problem, define control goals, and seek for optimal control actions, MPC and RL are representatives of model-based and learning-based control approaches, respectively. Although researchers have applied certain forms of MPC and RL to control the greenhouse climate, very few effort has been allocated to analyze connections, differences, pros and cons between MPC and RL either from a mathematical or performance perspective. Therefore, this paper will 1) propose MPC and RL approaches for greenhouse climate control in an unified framework; 2) analyze connections and differences between MPC and RL from a mathematical perspective; 3) compare performance of MPC and RL in a simulation study and afterwards present and interpret comparative results into insights for the application of the different control approaches in different scenarios. 

\end{abstract}


\end{frontmatter}


\clearpage
\section{Introduction}
\label{sec:introduction}
The world population has grown drastically in recent decades. Although its growth is slowing down, it is estimated that the world population will increase by 2,000 million people in the next 30 years reaching 9,700 million people in 2050~\cite{UN2019}. To feed this population in 2050, projections show that food production would need to increase by 70\% between 2005/07 and 2050. Production in developing countries would have to nearly double~\cite{FAO2018}. On the other hand, the Intergovernmental Panel on Climate Change (IPCC) declared that it is crucial and urgent to render the way of land use and the agricultural production methods more efficient in order to curb global warming. 

All of these facts have motivated innovative production methods and technical solutions in the agricultural sector to improve agri-food production and increase yield per hectare. Climate controlled greenhouses, which allows for growing crops regardless of the outdoor environment, is one of the important growing methods, especially in a changing climate. To maintain proper growing climate with efficient energy usage and operational cost, advanced control methods of greenhouse production system (e.g. lighting, heating, CO$_2$ dosing, ventilation, screening, etc.) are also needed. Moreover, as the number of greenhouse production systems is increasing, while the number of experienced growers is limited, autonomous climate control of a greenhouse production system is also necessary.

Among various control methods, MPC is effective to optimize a greenhouse production system due to its promising performance in multi-input and multi-output systems. However, predicting disturbances along the prediction horizon remins a challenging task~\cite{Boersma2021,Maciejowski2002,Chen2020}. The first implementations of MPC for greenhouse climate control can be traced back to the beginning of this century. The authors in~\cite{ElGhoumari2005} illustrate that a real-time applied MPC outperforms an adaptive PID controller, demonstrating its potential. MPC applications can also be found in, among others,~\cite{vanHenten1994,vanStraten2011,Pinon2005,Coulho2005,Blasco2007,Gruber2011,Chen2018,Kuijpers2022,Boersma2022}. Furthermore, \cite{Blasco2007} firstly introduces a non-linear MPC incorporating energy and water consumption to maintain climatic conditions in a greenhouse near the coast of Spain. Comparison results show that MPC can work better than set-point tracking controller, in some areas achieving an improvement of up to 10 times. \cite{Gruber2011} presents another nonlinear MPC approach based on a Volterra model \cite{Volterra1959} that captures the nonlinear relationship between ventilation and temperature. The proposed nonlinear MPC is finally applied to a detailed simulation greenhouse model. \cite{Gruber2014} also used a non-linear MPC to control the greenhouse temperature and its ventilation, as well as a hybrid MPC in \cite{Montoya2016}. Besides nonlinear MPC, linear MPC is also designed for temperature control of a greenhouse by \cite{Ioffe2015}, which behaves much better than the conventional on-off pulse-width modulation controller that was previously implemented. Moreover, \cite{Gonzalez2014} presents a tube-based linear MPC for a greenhouse system with two-time-scale dynamics. With development of sensing technology and data science, MPC tends to enter into a new stage where data-based techniques are integrated into control design. For instance, \cite{Chen2020} proposes a data-driven MPC for greenhouse climate control, mainly focusing on temperature and carbon dioxide concentration level. The main contribution of \cite{Chen2020} is to combine a dynamical model with a data-based model in order to identify uncertainties in the weather forecast. Besides optimization, \cite{Lin2020} also involved tracking performance using a hierarchical MPC improving energy efficiency and reduction of operational costs. More precisely, a two layer-based MPC structure was defined, where the upper layer generates optimal set-points for the greenhouse climate control. The lower layer is introduced to track the trajectories produced by the upper layer.


Recent advances in Information and Communication Technologies (ICT), as well as artificial intelligence (AI) have motivated the usage of AI-based or learning-based control, more specifically Reinforcement Learning (RL), as the key technology to transform the modern farming control mindset. Actually, the development and application of RL has a long and well-known history (see \cite{sutton1998} and references therein). In the last two decades, several applications have been published to use RL to solve real world challenges \cite{busoniu2010}, where the most relevant come from the adoption of the deep learning paradigm \cite{mousavi2016}. 
RL comprises a collection of algorithms and techniques that learn to solve a control problem by trial and error interaction with the environment (system). 
Due to this characteristic, the control model and strategy of the climate in a greenhouse could be updated and adaptable for different stages of the plant's development, even for difference cultivars. Therefore, RL is becoming more and more popular as an option for optimal and autonomous greenhouse climate control. A sign of the popularity of RL in greenhouse production systems is the three editions of the Autonomous Greenhouse Challenge \cite{Hemming2019,Hemming2020}, where dozens of international teams tested their state-of-the-art AI algorithms in greenhouse climate, irrigation, and crop growth control. Besides, greenhouse benchmarks and the so-called gyms (which are software environments to train machine learning agents) have also been published recently \cite{Overweg2021,Turchetta2022,An2021}. 
Among which, the first known reference in literature is \cite{tchamitchian2005daily} that describes a RL temperature controller in a rose production greenhouse. Other variables are not controlled and few details are given about the creation of the controller. In \cite{ban2017control}, the authors have already used a deep learning approach observing seven input variables from the environment and producing eight control actions (some of them are Boolean). This work is based on a greenhouse simulator and aims at keeping the variables within pre-established bounds. A more recent article is \cite{zhang2021robust}, where the aim is not only to control the system but also to create optimized models of the greenhouse dynamics using samples from the real environment. The focus of \cite{zhang2021robust} is on the robustness of the models and the controller. In \cite{afzali2021optimal}, a supplemental lightning control system is described, which is developed with Q-learning, posing the problem as a discrete constrained optimal problem where energy is the variable to optimize.

From the literature presented above, it is clear that certain forms of MPC and RL have already been applied to control greenhouse's climate. However, up to now, very few effort has been allocated to analyze connections, differences, pros and cons between MPC and RL especially for greenhouse climate control. However, it is important to do so as the comparisons may bring deep knowledge about these two interesting methods and it also brings insight on how to use and when to choose both methods in a greenhouse climate control application. 

Consequently, the objective of this paper is to compare the development and application of MPC and RL on greenhouse climate control. An agent-based deep RL controller is developed for a lettuce greenhouse system in a unified framework with a MPC approach, using the Deep Deterministic Policy Gradient (DDPG) approach. The used MPC approach to compare is mainly based on a nonlinear model, with multiple input and output variables, together with climate disturbances \cite{Boersma2021}. 
The outline of this paper is as follows: a lettuce greenhouse model is presented in Section~\ref{sec:model} on which the MPC and RL controllers are based. The proposed lettuce greenhouse model is also used for simulations and control strategy evaluation. Section~\ref{sec:controllers} details the development of the MPC and RL controllers from a unified framework. The simulation results based on these two different control approaches are comparatively presented in Section~\ref{sec:results}. Finally, Section~\ref{sec:conclude} discusses and concludes this work with insights for selecting and applying MPC and RL in different scenarios.

\section{Lettuce Greenhouse Model}
\label{sec:model}

\begin{figure}[!ht]
	\centering
	\includegraphics[width=0.8\textwidth]{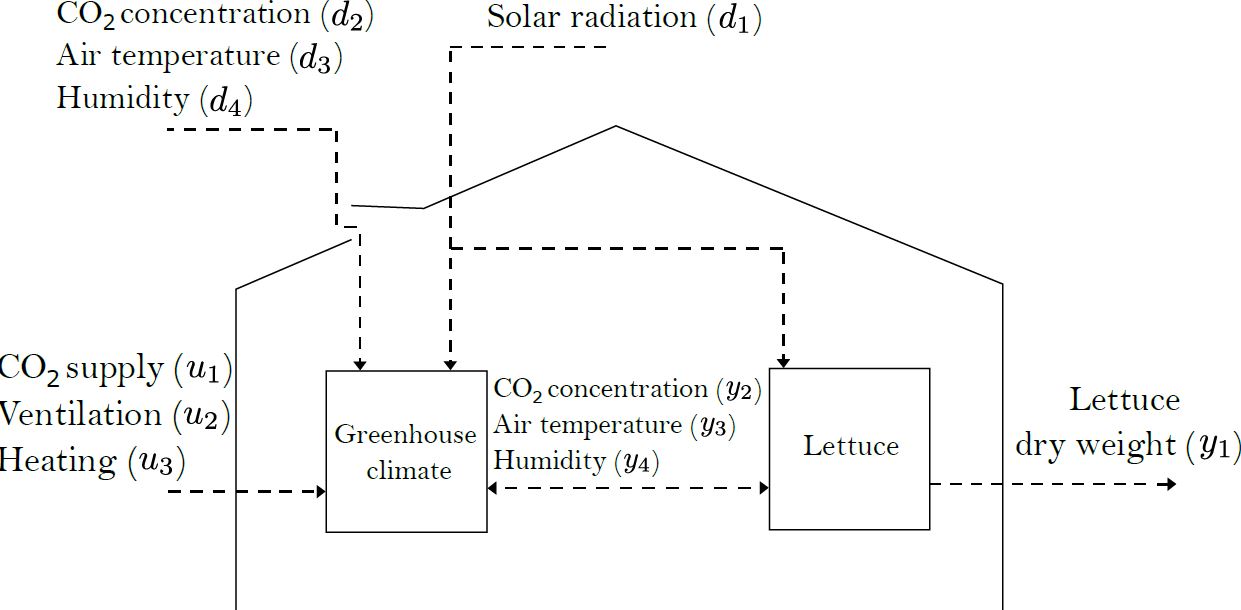}
	\caption{Schematic representation of lettuce greenhouse.} \label{fig:greenhouse}
\end{figure}

The lettuces greenhouse model is taken from~\cite{vanHenten1994} and discretized using explicit fourth order Runge-Kutta method with sample period $h$. Consequently, the following state-space model can be defined:
\begin{equation}
	\begin{aligned}
		x(k+1) &= f\big(x(k),u(k),d(k),p\big), \\
		y(k)   &= g\big(x(k),p\big), 
		\label{eq:model}
	\end{aligned}
\end{equation}
with time $k \in \mathbb{Z}^{0+}$, state $x(k) \in \mathbb{R}^4$, measurement $y(k) \in \mathbb{R}^4$, control input $u(k) \in \mathbb{R}^3$ and weather disturbance $d(k) \in \mathbb{R}^4$. Parameter $p \in \mathbb{R}^{28}$ and nonlinear functions $f(\cdot) g(\cdot)$ are given in Appendix. The state $x(k)$ contains the dry matter content of the lettuce $x_1(k)$ in kg$\cdot$m$^{-2}$, which is the lettuce's weight per square meter after all water has been removed. The state additionally contains the indoor $CO_2$ concentration $x_2(k)$ in kg$\cdot$m$^{-3}$, air temperature $x_3(k)$ in °C and humidity in $x_4(k)$ in kg$\cdot$m$^{-3}$. The weather disturbance $d(k)$ contains the incoming radiation $d_1(k)$ in W$\cdot$m$^{-2}$ and the outside $CO_2$ concentration $d_2(k)$ in kg$\cdot$m$^{-3}$, temperature $d_3(k)$ in °C and humidity content $d_4(k)$ in kg$\cdot$m$^{-3}$. The control signal $u(k)$ contains supply rate of $CO_2$ $u_{1}$ in mg$\cdot$m$^{-2}\cdot$s$^{-1}$, ventilation rate through the vents $u_2(k)$ in mm$\cdot$s$^{-1}$ and energy supply by heating system $u_3(k)$ in W$\cdot$m$^{-2}$. The measured output vector $y(k) \in \mathbb{R}^4$ contains $x_1(k)$ in g$\cdot$m$^{-2}$, $x_2(k)$ in ppm, $x_3(k)$ in °C and $x_4(k)$ in \%. Figure~\ref{fig:greenhouse} depicts the greenhouse model with lettuce~\eqref{eq:model}. 

\section{MPC versus RL}
\label{sec:controllers}
Model predictive control and reinforcement learning have been developed by separate communities, the control system community and the computational intelligence community, respectively \cite{Gorges2017}. Afterwards, the two methods have evolved more or less independently. As representatives of model-based and learning-based control methods, MPC and RL behave differently from the use of terminologies to the way of seeking for optimal control actions. 

Due to the way of formulating the control problem and defining the control goal, the implementation of a MPC requires a good model. The adaptability of MPC to various conditions is limited (assuming the model in the MPC is not online updated) and may not be enough for autonomous objectives in greenhouse production. Moreover, handling uncertainties in a MPC is computationally expensive due to the complex mathematical propagation of these uncertainties. 

Reinforcement learning is a dynamic control strategy which can automatically update the current control policy through incorporating newly developed knowledge learned from historical and real time data. Due to this characteristic, the control strategy of RL is more adaptable to current dynamics of the systems. Specifically for a greenhouse it is more adaptable to different stages of plant development, even for difference cultivars. Besides, RL can learn control strategies completely from data and not always a good model is needed. Another advantage is that it is relatively easy applied with limited complexity in manipulating the instruments. In spite of being simple and well-grounded, RL suffers from clear limitations. The most significant is the necessity to have a discrete and reduced set of problem states and control actions.  

 Apart from differences, MPC and RL also share plenty of common features. For example, both MPC and RL are predictive controllers independently of whether they integrate disturbance forecasts in their control logic. MPC uses explicit optimization along a finite prediction horizon, while RL learns actions to optimize the sum of the immediate and the discounted future rewards. In order to compare MPC with RL, the following subsections will explain precisely the development of RL from an unified framework with MPC in terms of 1) prediction principle; 2) reward (or cost) function; as well as 3) constraints handling. Besides the mathematical comparisons of the prediction principle, reward function and constraints handling, the performance of the developed MPC and RL, that are applied on the same use case, will also be compared. The theoretical comparisons about stability, feasibility, and robustness of the two control approaches will not be covered in this paper.

In the following, first the MPC framework is detailed (Section~\ref{sec:controller}) and then the RL framework is detailed (Section~\ref{sec:DDPG}).

\subsection{Nonlinear Model Predictive Control}
\label{sec:controller}

Based on the lettuce greenhouse model presented in Section~\ref{sec:model}, this subsection presents a nonlinear MPC for climate control to optimize the greenhouse's efficiency.

\subsubsection{Prediction Principle}

The proposed MPC controller is working according to the receding horizon principle. Here, at each time step, a new state measurement is taken from the greenhouse and used to initialize the model in the MPC (controller model). This model is propagated forward in time from this current state measurement, while a constrained cost function is minimized with the control inputs as decisions variables. From these optimized future control input sequences, only the first value is applied after which, again, a new current state measurement is collected. This procedure is repeated at each time step. The model in~\eqref{eq:model} is used in the MPC and in the following, the optimization problem, cost and constraints are formulated.

\subsubsection{Optimization Problem}

The optimization problem that is formulated in this section is employing the model given in~\eqref{eq:model}. It is assumed that at each time instant, the state $x(k)$ can be measured or is perfectly estimated. Then, the following optimization problem is solved at each time step $k_0$:
\begin{equation}
	\begin{aligned}
		\min_{u(k)} ~ & \sum_{k=k_0}^{k_0+N_p} V\big(u(k),y(k)\big), \\
		\text{s.t.} 	\quad & x(k+1) = f\big(x(k),u(k),d(k),p\big), \quad y(k) = g\big(x(k),p\big),\\
								~ & u_{\text{min}} \leq u(k) \leq u_{\text{max}}, \quad |u(k)-u(k-1)| \leq \delta u, \\
								~ & y_{\text{min}}(k) \leq y(k) \leq y_{\text{max}}(k), \quad \text{for}~k=k_0,\ldots,k_0+N_p,\\
								~ & x(k_0) = x_0.
	\end{aligned}
	\label{eq:optimization}
\end{equation}

\subsubsection{Cost Function and Constraints}

The cost function $V\big(u(k),y(k)\big)$ is defined as:
\begin{equation}
	\begin{aligned}  
		V\big(u(k),y(k)\big) 	&= -q_{y_1} \cdot y_1(k_0+N_p)+ \sum_{j=1}^{3} q_{u_j} \cdot u_j(k), 	 
	\end{aligned}
    \label{eq:MPC_cost_function}
\end{equation}
with $q_{y_1}, q_{u_j} \in \mathbb{R}$ defined as weights in the optimization that can also be seen as tuning variables. This cost function establishes a trade-off between the maximization of yield per square meter and the minimization of energy usage (control inputs). This trade-off is determined by the ratio of $q_{y_1}$ and $q_{u_j}$.

\begin{figure}[!ht]
	\centering
	\includegraphics[width=.6\textwidth]{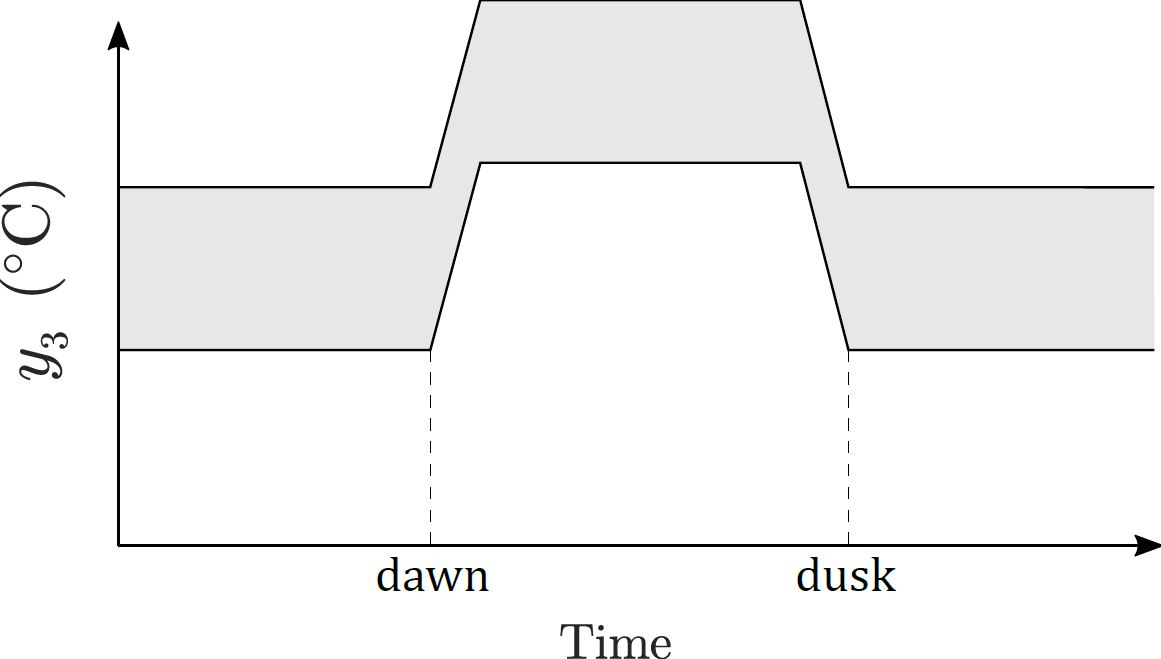}
	\caption{Graphical illustration of time-varying constraint imposed on output $y_3(k)$ that represents the temperature inside of the greenhouse. The gray area is the region where the optimized control signals steer $y_3(k)$ in.  \label{fig:TemperatureConstraint}}
\end{figure}

The constraints in~\eqref{eq:optimization} are defined as:
\begin{equation}
	\begin{aligned}  
		u_{\text{min}} &= \begin{pmatrix} 0 & 0 & 0 \end{pmatrix}^T, \quad u_{\text{max}} = \begin{pmatrix} 1.2 & 7.5  & 150 \end{pmatrix}^T, \quad \delta u = \frac{1}{10} u_{\text{max}}, \\ 
		y_{\text{min}}(k) &= \begin{pmatrix} 0 & 0 & f_{y_{3,\text{min}}}(k) & 0 \end{pmatrix}^T, \\ 
		y_{\text{max}}(k) &= \begin{pmatrix} \infty & 1.6 & f_{y_{3,\text{max}}}(k) &  70 \end{pmatrix}^T,
	\end{aligned}
\end{equation}
with lower and upper bounds on the control input defined by $u_{\text{min}}, u_{\text{max}} \in \mathbb{R}^3$, respectively, and the bound on the change of the control input defined by $\delta u \in \mathbb{R}^3$. The time-varying lower and upper bound on the output are $y_{\text{min}}(k)$ and $y_{\text{max}}(k) \in \mathbb{R}^4$, respectively. More precisely, only the third element in each of these bounds is time-varying and defined as:	
\begin{equation}
	\begin{aligned} 
		f_{y_{3,\text{min}}}(k) &= \begin{cases} 10, & \mbox{if } d_1(k_0) < 10 \\ 15, & \mbox{otherwise} \end{cases}, \quad	f_{y_{3,\text{max}}}(k) = \begin{cases} 15, & \mbox{if } d_1(k_0) < 10 \\ 20, & \mbox{otherwise.} \end{cases}
	\end{aligned}
	\label{eq:TemperatureConstraint}
\end{equation}
These time-varying constraints on the indoor temperature are set such that the indoor temperature is colder during the night than during the day in the greenhouse according to~\cite{Seginer1994}. Here it is demonstrated that lower greenhouse temperatures can later be compensated by higher ones as long as a daily average greenhouse temperature is satisfied. The time-varying constraint on the indoor temperature is graphically illustrated in Fig.~\ref{fig:TemperatureConstraint}. Here, the gray area indicates the region where the controller model output $y_3(k)$ is controlled to.   

\subsection{Deep Reinforcement Learning}
\label{sec:DDPG}
Deep Deterministic Policy Gradient (DDPG) is used to develop the RL control agent. The DDPG algorithm stems from Deterministic Policy Gradient \cite{Silver2014} algorithm comprising concepts of deep learning theory. The main advantage of DDPG is that it provides good performance in large and continuous state-action space environments, which motivates the selection in the greenhouse climate control application at hand.

DDPG is an off-policy and model-free actor-critic RL algorithm \cite{Lillicrap2016}. Due to these characteristics, the control actions of DDPG are generated by a different policy than the one being learnt and the optimal policy and value function are estimated directly without making efforts to learn the system dynamics.

\begin{figure}[!ht]
\centering 
\includegraphics[width=0.6\columnwidth]{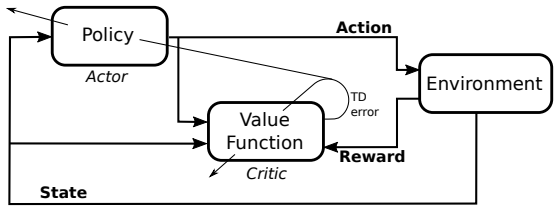}
\caption{Actor-Critic agent structure.}
\label{Fig:A-Cstruct}
\end{figure}

The structure of an actor-critic RL agent is shown in Figure \ref{Fig:A-Cstruct}, where the actor has the role of storing and applying the current best policy, using a deep neural network. According to the learned policy function, $\pi(s)$, and the system state, $s$, the actor computes the optimal actions, $u$. On another hand, the critic has the role of storing the value function, $Q(s,u)$, also using a deep neural network. The value function is the expected accumulated future reward for each state-action pair. The critic estimates the value function using the reward obtained from the system (environment) and its own information. Moreover, the critic is also in charge of calculating the temporal-difference error (TD) (i.e. the loss function), which is used during the learning process for both the critic and the actor. 

More than the actor and critic roles, DDPG also uses two distinctive elements of Deep-Q-Network \cite{Mnih2015}: the replay buffer and target networks. The replay buffer is a memory buffer that stores the transition tuple at each step. This tuple contains the current state $s(k)$, action $u(k)$, obtained reward $r(k)$, next state $s(k+1)$ and a Boolean variable indicating if the next state is terminal $t(k)$. A terminal state designates a state where the experiment ends. At each timestep, the critic and the actor are trained with a minibatch obtained by sampling random tuples from the replay buffer. This way of training eliminates time correlation between learning samples and facilitates convergence during the learning process.

Target networks are copies of the actor and critic networks. They are used during the training phase, providing the target values to compute the loss function. Once the original networks are trained with the set of tuples of the minibatch, the trained networks are copied to the target networks using a soft update, \textit{i.e.} forcing the target weights to change slowly. The use of target networks with soft update allows to give consistent targets during the TD backups and makes the learning process remain stable. Consequently, \textit{DDPG} requires four neural networks: the actor and the critic and their respective target networks.

The gradient functions that are used to update the weights of the critic and actor are presented in \eqref{Eq:update_Q} and \eqref{Eq:update_P}, respectively.  $\phi$ is the set of weights of the critic network and $\theta$ the weight of the actor, $\eta_{\phi}$ and $\eta_{\theta}$ are the learning rates of critic and actor, $B$ represents the mini-batch of transition tuples and $N$ is its size. Target networks are represented with the prime symbol. $\upsilon(k)$ \eqref{Eq:targets} are the target Q-values (not to be confused with target networks) and are used to compute the loss function. The weights of the critic are updated to minimize this loss function. The discount factor, $\gamma$, is a value between 0 and 1 that determines the importance of future rewards to the current state. Even though a one-to-one relation with the prediction horizon in MPC is difficult to define, there is a relation between these two. Note that the target Q-Values \eqref{Eq:targets} are obtained from the outputs of the actor and critic target networks, following the target network concept

\begin{equation}
\Delta \phi = \eta_{\phi}\nabla_{\phi} \left( \frac{1}{N} \sum_{k \in B} \left( Q(s(k),u(k)\mid \phi^{Q'})-\upsilon(k) \right)^2 \right)
\label{Eq:update_Q}
\end{equation}

\begin{equation}
\Delta \theta = \eta_{\theta}\nabla_{\theta} \left( \frac{1}{N} \sum_{k \in B}  Q(s(k),\pi(s(k)\mid \theta^{\pi})\mid \phi^{Q}) \right)
\label{Eq:update_P}
\end{equation}

\begin{equation}
\upsilon(k) = r(k) + \gamma Q'(s(k+1),\pi'(s(k+1) \mid \theta^{\pi'})\mid \phi^{Q'})
\label{Eq:targets}
\end{equation}

The update of the weights of the target networks from the trained networks are shown in \eqref{Eq:soft_critic} and \eqref{Eq:soft_actor}. The parameter $\tau$ indicates how fast this update is carried on. This soft update is made each step after training the main networks

\begin{equation}
\phi^{Q'} \leftarrow \tau \phi^{Q} + (1- \tau) \phi^{Q'}.
\label{Eq:soft_critic}
\end{equation}

\begin{equation}
\theta^{\pi'} \leftarrow \tau \theta^{\pi} + (1- \tau) \theta^{\pi'}.
\label{Eq:soft_actor}
\end{equation}

\subsubsection{Agent Description}
The agent's behaviour is shaped according to it's state, the reward function and the discount factor. The state contains ten terms:
\begin{multline}
s(k)=(\Delta_{y_1}(k), e_{y_2}(k), e_{y_3}(k), y_4(k), d_1(k), d_2(k), d_3(k), \\ u_1(k-1), u_2(k-1), u_3(k-1))
\label{Eq:agent_state}
\end{multline}
where $\Delta_{y_1}(k)=y_1(k)-y_1(k-1)$, $e_{y_2}(k)=y_{2,ref}(k)-y_2(k)$, $e_{y_3}(k)=y_{3,ref}(k)-y_3(k)$, $d_i(k)$ are the current disturbances and $u_i(k-1)$ are the previous control actions.
The reward function and the discount factor take similar roles as the cost function and prediction horizon, respectively, in MPC. They are detailed next.

\paragraph{\textbf{Reward Function and Constraints}}
Similarly to the cost function of MPC presented in equation \eqref{eq:MPC_cost_function}, the dry matter is the most important item in the reward function. Nevertheless, other variables also need to be taken into account to guide the learning process. As a result, the reward function contains six terms. The first three terms are rewards and the other terms are penalties. It is defined as: 
\begin{equation}
r(k)=c_{r,1}\Delta_{y_{1}}(k) + r_{CO_2}(k)+r_{T}(k)
-(\sum_{j=1}^{3} c_{r,u_j} \cdot u_j(k-1))
\label{Eq:reward}
\end{equation}
where $r_{CO_2}(k)$ and $r_{T}(k)$ represent the rewards associated with the control of the main production variables, $CO_2$ concentration and temperature, and $c_{r,i}$ are adjustable constant parameters.

Dry matter and energy consumption are considered as in the MPC formula \eqref{eq:MPC_cost_function} but the agent is also rewarded if $CO_2$ concentration and temperature are controlled satisfactorily. In accordance with the temperature constraint included in the MPC \eqref{eq:TemperatureConstraint}, the reward function depends not only on the temperature error but also on the heating system. That is because during the day it is more efficient to warm up the greenhouse with the solar radiation than with the heating system and applying ventilation. $CO_2$ concentration is also controlled to be higher during the day and lower during the night, as is common practice in a greenhouse.

The agent will receive a fixed reward for maintaining these measurements within a given range and a penalty for going outside of the range. Consequently, the constraints in \eqref{eq:optimization} are codified here as rewards, defined as

\begin{equation}
\begin{aligned} 
r_{CO_2}(k) &= \begin{cases} 
-c_{r,CO_2,1} \cdot (y_2(k)-CO_{2_{min}}(k))^2 & \mbox{if  } y_2(k) < CO_{2_{min}}(k)  
 \\ -c_{r,CO_2,1} \cdot (y_2(k)-CO_{2_{max}}(k))^2 & \mbox{if  } y_2(k) > CO_{2_{max}}(k)
 \\ c_{r,CO_2,2} & \mbox{otherwise}
\end{cases} 
\cr
 r_{T}(k) &= \begin{cases} 
-c_{r,T,1} \cdot (y_3(k)-T_{min}(k))^2 & \mbox{if  } y_3(k) < T_{min}(k)  
 \\ -c_{r,T,1} \cdot (y_3(k)-T_{max}(k))^2 & \mbox{if  } y_3(k) > T_{max}(k)
 \\ c_{r,T,2} & \mbox{otherwise},
\end{cases} 
\cr
\end{aligned}
\label{Eq:sub_reward}
\end{equation}
with $c_{r,CO_2}$ and $c_{r,T}$ being adjustable constant parameters. These parameters can be considered as the weights of the constraint on the reward. Depending on the strength of the reward, the agent will adopt a more conservative or more aggressive strategy for the greenhouse climate control.

To sum up, the reward function of the RL agent codifies the MPC cost function and constraints. In \eqref{Eq:reward}, the dry matter increase is rewarded instead of the accumulated value, as in \eqref{eq:MPC_cost_function}, because the agent learning was not consistent when the accumulated value was used. 

\paragraph{\textbf{Prediction Principle}}
Not like the receding horizon used in MPC, the RL agent uses a discount factor $\gamma$, as shown in \eqref{Eq:targets}, to describe the future influence of rewards. If $\gamma = 0$ , the RL agent will completely focus on learning optimal actions for the immediate reward, while if $\gamma = 1$, the RL agent will evaluate its actions based on the total sum of all its future rewards. So that picking a particular value of $\gamma$ is equivalent to picking a prediction horizon of MPC though it is complicated to find an exact one-to-one relation.

 
 From \eqref{Eq:targets}, one can also see that discount factor $\gamma$ is related to the aggressiveness of the control actions, similarly to parameters $q_{y_1}$ and $q_{u_j}$ in \eqref{eq:MPC_cost_function}. In this paper, the discount factor $\gamma$ is tuned by trial and error to produce dynamical behaviour in the greenhouse similar to MPC.
 
\subsubsection{Agent training} 

\begin{figure}[h!]
\includegraphics[width=1.0\textwidth]{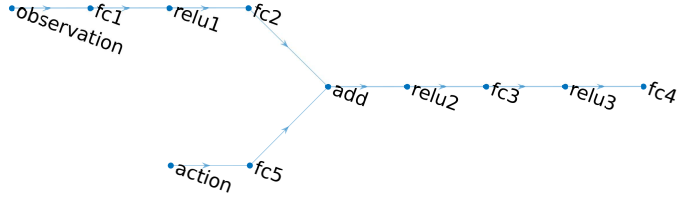}
\caption{Actor-Critic agent structure.}
\label{Fig:criticStructure}
\end{figure}

The structure of the RL agent developed in this paper is described in this section. The critic networks structure is shown in Figure~\ref{Fig:criticStructure}. Observations are processed by three layers of 10 Rectified Linear Units (ReLU) and actions are processed by two layers of 10 units. All the layers are fully connected. The actor networks consist of three fully connected layers of 20 ReLU units and a final hyperbolic tangent layer with three nodes, corresponding to each one of the control actions. The training options for all the networks are given in Table~\ref{tab:training_options} and the parameters related to the reward/penalty function are listed in Table~\ref{tab:agent_parameters}.

\begin{table}[h]
	\center
	\caption{\label{tab:training_options}Training options for all the networks}
	\resizebox{0.85\textwidth}{!}{	\begin{tabular}{ cc | cc} 
			parameter & value & parameter & value\\
			\hline
			learning rate & \num{1e-3} & experience buffer size 	& \num{1e4}\\ 
			gradient threshold & \num{1}  & experience mini-batch size & 64\\ 
			$L_2$ regularization factor & \num{1e-5} & discount factor, $\gamma$ 	& 0.9 \\ 
	\end{tabular}}
\end{table}	

The agent is trained for 500 epochs, where each epoch consists of one day of crop growth simulation. Each day of simulation contains 96 steps, so each step is 15 minutes. The agent training uses different meteorological conditions and initial conditions for each epoch to avoid overfitting and to facilitate generalization. Therefore disturbances are $d_i(k)=d_i(k)\cdot \kappa$ where $\kappa \sim U(0.7, 1.3)\,$.

\begin{table}[h]
	\center
	\caption{\label{tab:agent_parameters}Parameters of the reinforcement function}
	\resizebox{0.6\textwidth}{!}{	\begin{tabular}{ cc|cc } 
			parameter & value & parameter & value\\
			\hline
			$c_{r,1}$ & 16 & $c_{r,CO_2,1}$ 	& 0.1\\ 
			$c_{r,u_1}$ & \num{-4.5360e-04}  & $c_{r,CO_2,2}$ & 0.0005\\ 
			$c_{r,u_2}$ & \num{-0.0075} & $c_{r,T,1}$ 	& 0.001 \\ 
			$c_{r,u_3}$ & \num{-8.5725e-04} & $c_{r,T,2}$	& 0.0005\\ 
	\end{tabular}}
\end{table}	

\section{Simulation Results}
\label{sec:results}
\begin{figure}[h!]
\includegraphics[width=1.0\textwidth]{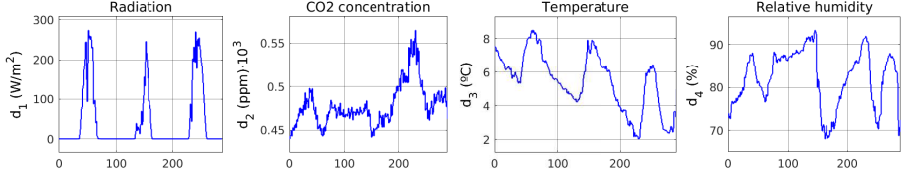}
\caption{Disturbances applied to the greenhouse (solar radiation, external $CO_2$ concentration, temperature and relative humidity).}
\label{Fig:disturbances}
\end{figure}

In order to be able to compare MPC with RL, the same disturbance is used and similar constraints and control goals are taken into account The weather data $d(k)$ used throughout the simulations are real-life data, presented in~\cite{Kempkes2014}. These data are collected during experiments performed in the greenhouse called ``the Venlow Energy greenhouse'' that is located in Bleiswijk, Holland. The collected data points are sampled at 5 minutes and $N$ of these are used and re-sampled to the sample period $h$. Figure~\ref{Fig:disturbances} shows the specific disturbance realization used to compare the performance of the DDPG agent and the MPC controller for 288 samples, which are equivalent to 3 days.
\subsection{Economic Profit Indicator} 
Besides dry matter, indoor $CO_2$ concentration, air temperature and relative humidity are optimized in the cost (reward) function, also one extra indicator is considered. This is the economic profit indicator (EPI) and considered as an additional validation item to compare the performance of the MPC and the DDPG-based RL agent:
\begin{equation}
EPI=\phi({y_1(t_f)}) -\sum_{t_b}^{t_f} (c_q u_q(t) + c_{co_2}u_{co_2}(t)) h,
\label{Eq:profit}
\end{equation}
\noindent where $\phi({y_1(t_f)})$ is the gross income obtained at harvest time $t_f$ and $c_q u_q(t) + c_{co_2}u_{co_2}(t)$ are the operating costs of the air conditioning equipment ($Hfl·m^{-2}·s^{-1}$). The auction price of the lettuce follows a linear ratio $\phi({y_1(t_f)}) = c_{pri,1} + c_{pri,2}y_1(t_f)$, between the auction price and the harvest weight of lettuce, in $kg·m^{-2}$. The units of parameters $c_{pri,1}$ and $c_{pri,2}$ are $Hfl·m^{-2}$ and $Hfl·kg^{-1}·m^{-2}$, respectively.

It is assumed that the operating costs of the climate control equipment are linearly related to the amount of energy $u_q$ ($W·m^{-2}$) and the amount of carbon dioxide introduced into the system is $u_c$ ($kg·m^{-2}·s^{-1}$). These operating costs are parameterized by the price of energy $c_q$ ($Hfl·J^{-1}$) and the price of carbon dioxide is $c_{co_2}$ ($Hfl·kg^{-1}$). More details of their values are found in Table~\ref{tab:profit}.

\begin{table}[h]
	\center
	\caption{\label{tab:profit}Parameters of the economic profit function}
	\resizebox{0.4\textwidth}{!}{	\begin{tabular}{ c|c } 
			parameter & value \\
			\hline
			$c_{co_2}$ & \num{42e-2}$Hfl kg^{-1}$\\ 
			$c_q$ & \num{6.35e-9}$Hfl J^{-1}$\\ 
			$c_{pri,1}$ & \num{1.8}$Hfl m^{-2}$\\ 
			$c_{pri,2}$ & \num{16}$Hfl kg^{-1}$\\ 
	\end{tabular}}
\end{table}	

\subsection{Specific MPC settings} 
 The weights $q_{\hat{y}_1},q_{u_i}$ are tuned such that an acceptable trade-off between yield and energy usage is achieved. The prediction horizon $N_p$ is not taken too large to prevent the necessity of including uncertainty that grows over time. Indeed, weather forecasts become more uncertain over the future horizon. Other settings that are used during the simulation studies are given in Table~\ref{tab:simulation_parameters}.
\begin{table}[h]
	\center
	\caption{\label{tab:simulation_parameters} Simulation and controller settings.}
    \resizebox{0.6\textwidth}{!}{	 
	\begin{tabular}{ cc|cc} 
		parameter & value  &parameter & value  \\
		\hline
		$h$ 									& 15 minutes  				& $q_{\hat{y}_1}$ 				& $10^3$			\\ 
		$N_p$ 									& 6 hours	 		 	& $q_{u_i}$  					& $\{10,1,1\}$	\\ 
		$N$                                     & 40 days      		 	& $N_s$                         &  20 \\
	\end{tabular}}
\end{table}

The open-source software CasADi~\cite{Andersson2019} and solver IPOPT~\cite{Wachter2006} are used in a Matlab environment to solve the optimization problem formulated in~\eqref{eq:optimization}, while following the direct single-shooting method and warm start option of IPOPT.

\subsection{Results}
\begin{figure}[h!]
\includegraphics[width=1.0\textwidth]{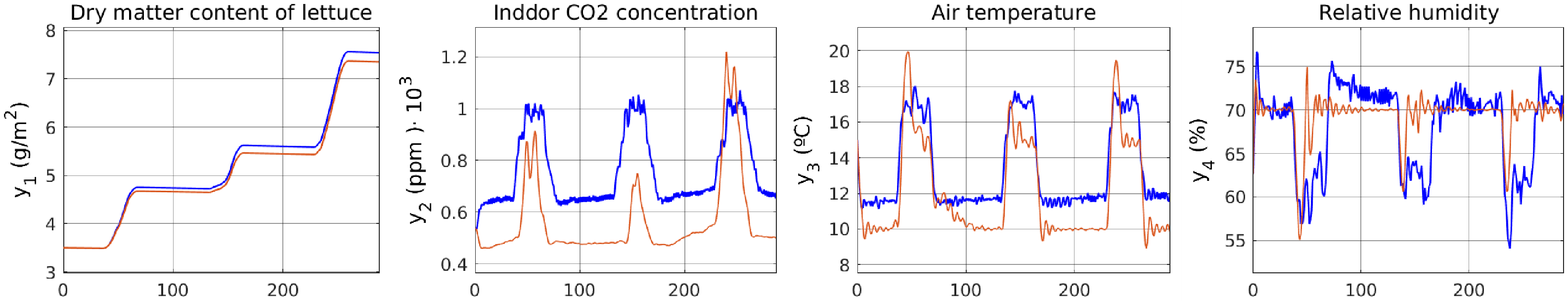}
\caption{Outputs of the system: dry matter, indoor $CO_2$ concentration, air temperature and relative humidity. DDPG agent (blue) and MPC (orange)}
\label{Fig:outputs}
\end{figure}

Figure~\ref{Fig:outputs} shows the measurement of the simulated greenhouse using the DDPG-based RL agent (in blue) and the MPC controller (in orange). From this figure, we can find that both MPC and RL obtain a similar dry matter content of lettuce, though RL is slightly more productive. The rest of the outputs are kept approximately within their constraints. The interior relative humidity levels are very similar in both cases, although RL allows for more variation. Similarly, the indoor temperature and $CO_2$ concentration are close to their minimum level although the RL agent keeps the temperature at a level slightly above the minimum temperature to avoid receiving this penalty. The RL agent is clearly more conservative than the MPC because the cost function penalizes when the temperature drops below the minimal temperature constraints. The MPC controller controls more accurately the humidity while the DDPG agent does a better job regarding the indoor temperature. The $CO_2$ concentration is increased noticeably during daylight hours by both controllers, which is expected. However, it is kept higher during night by the DDPG agent.

\begin{figure}[h!]
\includegraphics[width=1.0\textwidth]{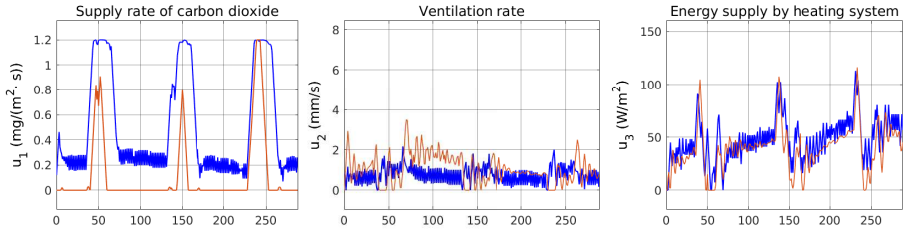}
\caption{Control actions applied by the DDPG agent (blue) and the MPC (orange).}
\label{Fig:controlactions}
\end{figure}

Figure~\ref{Fig:controlactions} shows the control actions of both controllers. It is apparent that more efficient use of ventilation is done by the MPC, but a higher quantity of energy is given by the DDPG agent. 

\begin{table}[h]
	\center
	\caption{\label{tab:perform}Performance of MPC and RL}
	\resizebox{0.6\textwidth}{!}{	\begin{tabular}{ c|c|c } 
			Index & MPC & RL \\
			\hline
			EPI & \num{1,843}$Hfl m^{-2}$ & \num{1,788}$Hfl m^{-2}$ \\ 
			computational time & \num{305, 04}$s$ &  \num{2,36}$s$ \\
	\end{tabular}}
\end{table}	


In terms of production, the RL agent achieves a greater production of lettuce, but in terms of economic benefit, the MPC achieves a higher economic return (\num{1,843} $Hfl m^{-2}$) than the agent (\num{1,788} $Hfl m^{-2}$), as shown in Table~\ref{tab:perform}. This is due to the fact that the economic profitability of the increase in lettuce production is lower than the cost of resource consumption for the economic return function with the parameters given in Table~\ref{tab:profit}. Figure~\ref{Fig:controlactions} shows that the temperature levels are kept low in both cases but the RL agent clearly keeps the but $CO_2$ concentration at a much higher level than the MPC. An explanation for this behaviour might be that since the increase in the injection of CO$_2$ can significantly increase the production, the agent has favored exploring this option arriving at a sub-optimal solution.

To solve this problem, we can let the agent train continuously to explore a larger state space, or we can reshape the cost function and redefine the upper and lower limits of the CO$_2$ injection. The agent can potentially find the optimal solution for this redefined problem. However, if the auction price of lettuce or the cost of CO$_2$ changes, this agent will not adapt to find a new optimal solution unless it takes all these variable parameters as observation and starts training again.
\begin{figure}[h!]
\includegraphics[width=1.0\textwidth]{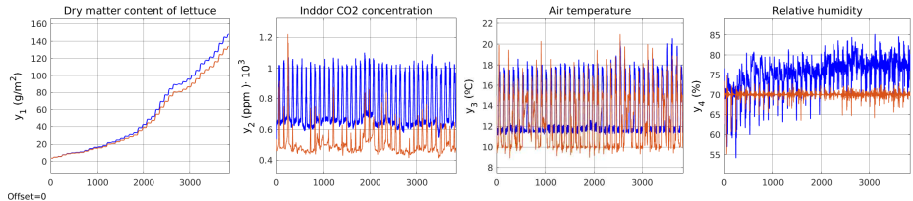}
\caption{System outputs for a complete growing cycle (40 days).}
\label{Fig:outputs40days}
\end{figure}

Figure~\ref{Fig:outputs40days} shows the outputs of the simulated greenhouse for a complete growing cycle of lettuce, which lasts for 40 days. The production of RL is now clearly higher, but the EPI is smaller (\num{2,195} $Hfl m^{-2}$) in contrast to \num{2,504} $Hfl m^{-2}$ produced by the MPC. It is also apparent that the RL agent is more permissive with the constraint of humidity. That might become a health problem for the crops and has to be looked at with care.


\section{Discussion and Conclusion}
\label{sec:conclude}
This paper proposed a MPC and RL-based control agent that control a lettuce greenhouse. The provided greenhouse model behaves as a simulation environment. The performance of the RL and MPC has been compared and analyzed in order to explore insights of using two different control methods for optimal control of greenhouse production to achieve the best use of natural resources and production efficiency in the presence of uncertainty in the forecast of the future climate. 

Regarding the results, the agent has not obtained as good results as the MPC model in terms of energy optimization. Furthermore, although more production has been obtained, the economic profitability has not been as good as in the case of the MPC. As in greenhouse, it is difficult to evaluate the controller in terms of output due to the large number of outputs, that is, it is difficult to design a cost function that can be accurately judged. Although we could directly use the economic profit function as a reinforcement function for the agent, this is not a reliable method for the following reasons.

In reinforcement learning we cannot directly apply constraints to the states, which leads to the fact that if we do not design a small and reasonable range for these states in the reward function, the agent will be able to explore those states that are completely irrational, such as going up the temperature above 40 ºC, turn on the fan to the maximum, etc. The second is that if there is a shortcut hidden in the cost function, the agent will go for that shortcut. For example, in the case of comparing the economic benefit, if the fan is set to maximum without injecting any carbon dioxide and without turning on the heating, the greatest economic benefit is obtained, which reaches \num{1,894} $Hfl m^{-2}$. However, this keeps the indoor temperature and relative humidity at a level that is very dangerous for the lettuces. Furthermore, with a trained agent, it is difficult to know whether the agent has converged to optimal control or suboptimal control. The only solution can be to let it compare itself with other optimal control controllers. In short, it is very difficult to design a well-tuned agent, and the design process is often an iterative one. 

However, once an agent that can be considered correct is obtained, it is very powerful, and control results can be obtained with the application of very little computing resources and computation time. Moreoever, RL agent can also handle uncertainties much easier than MPC in terms of its optimization strategies. To give clear comparisons between MPC and RL, we list all their cons and pros in the following Table~\ref{tab:compare}, which will give more insights on how to choose these methods for different scenarios.

\begin{table}[h]
	\center
	\caption{\label{tab:compare}Pros and Cons of MPC and RL}
	\resizebox{0.8\textwidth}{!}{	\begin{tabular}{ c|c } 
			 MPC & RL \\
			\hline
			Rely on good model & Can be model free \\
   Optimize from a quadratic convex model & Learning for decision making\\ 
   		Handle constraints easily & Difficulty to handle constraints \\ 
			Struggle with long term prediction & Infinite prediction horizons \\ 
   	Needs high computation load for uncertainties & With inherent robustness \\
       	Low adaptability & High adaptability \\
    	Online optimization complexity is high & Low complexity for online learning \\
	\end{tabular}}
\end{table}	

During the development of the project, we have identified a series of areas for improvement and continuation of the work:

\begin{enumerate}
     \item One of the proposals for continuing this work would be to redesign the agent and integrate future climate predictions as observations. In this work, we only input the current outdoor weather to the agent, and the prediction of the agent's future states is based solely on the current indoor and outdoor weather and the actions performed. In future work, one could try to include predictions of the future outdoor weather as input, thus allowing the agent to make better predictions of future reinforcements;
	\item  Another proposal is to use the parameters of the economic profit indicator as agent observations and let the agent learning to make corresponding changes in energy use when these parameters change;
	\item It is also proposed to apply this problem to more complex reinforcement deep learning algorithms, such as TD3 (Twin Delayed Deep Deterministic policy gradient algorithm), PPO (Proximal Policy Optimization) or SAC (Soft Actor-Critic), and compare their performance.
    \item Regarding both MPC and RL have their own cons and pros, another way to use both of their strong features is to integrate MPC with RL for climate control in greenhouse production system.
\end{enumerate} 

\clearpage
\section*{APPENDIX}
\noindent The greenhouse with lettuce model is defined as:
\begin{equation*}
	\begin{aligned}
		\frac{\text{d} x_1(t)}{\text{d}t} &=  p_{1,1} \phi_{\text{phot,c}}(t) -  p_{1,2} x_1(t) 2^{x_3(t)/10-5/2},  \\
		\frac{\text{d} x_2(t)}{\text{d}t} &= \frac{1}{p_{2,1}}  -\phi_{\text{phot,c}}(t) + p_{2,2} x_1(t) 2^{x_3(t)/10-5/2} + u_1(t) 10^{-6} - \phi_{\text{vent,c}}(t), \\
		\frac{\text{d} x_3(t)}{\text{d}t} &= \frac{1}{p_{3,1}} u_3(t) - (p_{3,2} u_2(t)10^{-3} +p_{3,3}) (x_3(t)-d_3(t)) + p_{3,4} d_1(t), \\
		\frac{\text{d} x_4(t)}{\text{d}t} &= \frac{1}{p_{4,1}} \big( \phi_{\text{transp,h}}(t) - \phi_{\text{vent,h}}(t) \big), \\
	\end{aligned}
\end{equation*}
with
\begin{equation*}
	\begin{aligned}
		\phi_{\text{phot,c}}(t) &= \Big(1-\text{exp} \big(-p_{1,3}x_1(t) \big) \Big) \Big(p_{1,4}d_1(t)\big(-p_{1,5}x_3(t)^2+ ... \\ 
		& \qquad p_{1,6} x_3(t) - p_{1,7} \big) \big( x_2(t)- p_{1,8}\big)\Big) /\varphi(t),\\
		\varphi(t)              &= p_{1,4}d_1(t)+\big(-p_{1,5}x_3(t)^2+p_{1,6} x_3(t) - p_{1,7} \big) \big( x_2(t)- p_{1,8}\big), \\
		\phi_{\text{vent,c}}(t) &= \big(u_2(t)10^{-3}+p_{2,3}\big)\big(x_2(t)-d_2(t)\big),\\
		\phi_{\text{vent,h}}(t) &= \big(u_2(t)10^{-3}+p_{2,3}\big)\big(x_4(t)-d_4(t)\big), \\
		\phi_{\text{transp,h}}(t) &= p_{4,2} \Big(1-\text{exp} \big(-p_{1,3}x_1(t) \big) \Big) \\ 
		& \qquad \Big( \frac{p_{4,3}}{p_{4,4}(x_3(t)+p_{4,5})}  \text{exp} \Big( \frac{p_{4,6} x_3(t)}{x_3(t) + p_{4,7}} \Big) - x_4(t) \Big),		
	\end{aligned}
\end{equation*}
and with $t \in \mathbb{R}$ the continuous time. Here, $\phi_{\text{phot,c}}(t), \phi_{\text{vent,c}}(t),$ $\phi_{\text{transp,h}}(t)$ and $\phi_{\text{vent,h}}(t)$ are the gross canopy photosynthesis rate, mass exchange of CO$_2$ through the vents, canopy transpiration and mass exchange of H$_2$O through the vents, respectively. The measurement equation is defined as:
\begin{equation*}
	\begin{aligned}
		y_1(t) &= 10^3 \cdot x_1(t)  	&\quad \text{g m$^{-2}$}, \\
		y_2(t) &= \frac{10^3 \cdot p_{2,4} \big( x_3(t) + p_{2,5} \big) } {p_{2,6} p_{2,7}} \cdot x_2(t),  	&\quad \text{ppm}  \cdot 10^3, \\
		y_3(t) &= x_3(t), 	&\quad \text{°C}, \\
		y_4(t) &=  \frac{10^2 \cdot p_{2,4} \big( x_3(t)+ p_{2,5} \big)}{11 \cdot \text{exp}\Big( \frac{p_{4,8} x_3(t)}{x_3(t)+p_{4,9}} \Big)} \cdot x_4(t), 	&\quad \text{\%},
	\end{aligned}
\end{equation*}
The model parameters $p_{i,j}$ are chosen following~\cite{vanHenten1994} and given in Table~\ref{tab:model_parameters}. 

\begin{table}[h]
	\center
	\caption{\label{tab:model_parameters}Values of the model parameters that are taken from~\cite{vanHenten1994}.}
	\resizebox{1\textwidth}{!}{	\begin{tabular}{ cc|cc|cc|cc } 
			parameter & value & parameter & value & parameter & value & parameter & value\\
			\hline
			$p_{1,1}$ & 0.544 & $p_{2,1}$ 	& 4.1 			& $p_{3,1}$ & 3$\cdot 10^{4}$ & $p_{4,1}$ 	& 4.1 \\ 
			$p_{1,2}$ & 2.65 $\cdot 10^{-7}$ & $p_{2,2}$ 	& 4.87 $\cdot 10^{-7}$ 			& $p_{3,2}$ & 1290 & $p_{4,2}$ & 0.0036 \\ 
			$p_{1,3}$ & 53 & $p_{2,3}$ 	& 7.5 $\cdot 10^{-6}$ 			& $p_{3,3}$ & 6.1 & $p_{4,3}$ 	& 9348 \\ 
			$p_{1,4}$ & 3.55 $\cdot 10^{-9}$ &	$p_{2,4}$	& 8.31  		& $p_{3,4}$ & 0.2 & $p_{4,4}$ 	& 8314 \\ 
			$p_{1,5}$ & 5.11 $\cdot 10^{-6}$ & $p_{2,5}$	& 273.15		&  			&   & $p_{4,5}$ 	& 273.15 \\ 
			$p_{1,6}$ & 2.3 $\cdot 10^{-4}$ & $p_{2,6}$	& 101325		&  			&   & $p_{4,6}$ 	& 17.4 \\ 
			$p_{1,7}$ & 6.29 $\cdot 10^{-4}$ & $p_{2,7}$	& 0.044			&  			&   & $p_{4,7}$ 	& 239 \\ 
			$p_{1,8}$ & 5.2 $\cdot 10^{-5}$ & 			&   			&  			&   &  $p_{4,8}$  	&   17.269  \\ 
			&   & 			&   			&  			&   							&  $p_{4,9}$	&  238.3  \\
	\end{tabular}}
\end{table}	

\noindent The model is discretized using the explicit fourth order Runge-Kutta method resulting in the discrete-time model as presented in~\eqref{eq:model}:
\begin{equation}
	\begin{aligned}
		x(k+1) &= f(x(k),u(k),d(k),p), \\
		y(k)   &= g(x(k),p),
	\end{aligned}	
\end{equation}
with discrete time $k \in \mathbb{Z}^{0+}$ and relation $t=k \cdot h$ with $h$ the sample period. The initial state and control signal that are used during the simulation are defined as: 
\begin{equation*}
	x(0) = \begin{pmatrix} 0.0035 & 0.001 & 15 &0.008 \end{pmatrix}^T, \qquad u(0) = \begin{pmatrix} 0 & 0 & 0 \end{pmatrix}^T. 
\end{equation*}

\clearpage
\bibliography{mybibfile}

\end{document}